\documentclass[12pt]{article}
 \usepackage[utf8]{inputenc}

      \usepackage{latexsym}
         \usepackage[reqno, namelimits, sumlimits]{amsmath}
         \usepackage{amssymb, amsfonts}
         \usepackage{amsthm}
           
 \newtheorem{theorem}{Theorem}[section]

 \newtheorem{remark}[theorem]{Remark}

 \newtheorem{pro}[theorem]{Proposition}

\author{ Gregory Seregin
}

\title{A note on  certain  scenarios of Type II blowups of suitable weak solutions to the Navier-Stokes equations  }

\author{G.~Seregin\footnote{University of Oxford, Mathematical Institute, OxPDE, Oxford, UK and St Petersburg Department of Steklov Mathematical Institute, RAS, Russia, email address: \texttt{seregin@maths.ox.ac.uk}}
}

\begin{document}


\maketitle


\begin{abstract} In the note, various scenarios of potential Type II blowups of suitable weak solutions to the Navier-Stokes equations are studied. It is shown, that under some assumptions, such type of blowups cannot happen. In this case, corresponding statements may be interpreted as regularity results. Their justification is based on a technique making use of a certain Euler scaling and Liouville type theorems for ancient solutions to the Euler system.
\end{abstract}


{\bf Data availability statement}
Data sharing not applicable to this article as no datasets were generated or analysed during the current study.

{\bf Acknowledgement} The work is supported by Leverhulme Emeritus Fellowship 2023. 


\setcounter{equation}{0}
\section{Introduction}

In this note, we are going to discuss different scenarios of Type II blowups of solutions to the Navier-Stokes equations on a local level. To this end, we need to introduce several notions and some notation. First, let us start with the definition 
of suitable weak solutions to those equations,  
see \cite{CKN}, \cite{Lin}, and \cite{LS1999}. A pair of functions $v$ and $q$ is called a suitable weak solution to the Navier-Stokes equations in space-time domain $Q$ if it has the following properties:
\begin{itemize}
	\item $v\in L_\infty(-1,0;L_2(B)),\quad \nabla v\in L_2(Q),\quad q\in L_\frac 32(Q);$
	\item the Navier-Stokes equations
$$\partial_tv+v\cdot\nabla v-\Delta v+\nabla q=0,\quad{\rm div}\,v=0$$
are satisfied in $Q$ in the sense of distributions;
\item
for a.a. $t\in ]-1,0[$, the local energy inequality 
$$\int\limits_B\varphi(x,t)|v(x,t)|^2dx+2\int\limits^t_{-1}\int\limits_B\varphi|\nabla v|^2dxd\tau\leq 
$$
$$\leq \int\limits^t_{-1}\int\limits_B (|v|^2(\partial_t\varphi+\Delta\varphi)+v\cdot\nabla\varphi(|v|^2+2q))dxd\tau$$ holds
for all smooth non-negative functions $\varphi$ vanishing in a vicinity of a parabolic boundary of the cylinder $Q$. 
\end{itemize}
Here and in what follows, the standard notation is used:
$Q(r)=B(r)\times ]-r^2,0[$ is a parabolic cylinder 
and $B(r)$ is a spatial ball of radius $r$ centred at the origin $x=0$, $B=B(1)$, and $Q=Q(1)$.

According to the partial regularity theory for suitable weak solutions 
 developed in the remarkable Caffarelli-Kohn-Nirenberg paper \cite{CKN}, there is a universal constant $\varepsilon_*>0$ such that if 
\begin{equation}
	\label{partial regularity}
	g(v)=\min\{\liminf_{r\to 0}A(v,r),\liminf_{r\to 0}E(v,r),\liminf_{r\to 0}C(v,r)\}<\varepsilon_*,
	\end{equation}
then $z=(x,t)=0$ is a regular point of $v$, i.e., $v\in L_\infty(Q(r_*))$ for a radius $r_*\in ]0,1[$. Here, the following quantities  have been used:
$$E(v,r)=\frac 1{r}\int\limits_{Q(r)}|\nabla v|^2dz, \quad A(v,r)=\sup\limits_{-r^2<t<0}\frac 1{r}\int\limits_{B(r)}|v(x,t)|^2dx,$$
$$C(v,r)=\frac 1{r^2}\int\limits_{Q(r)}|v|^3dz.$$
It remains unknown  whether or not $z=0$ is a regular point of $v$ for any suitable weak solution $v$ and $q$ in $Q$. 

Assuming that there exists a suitable weak solution $v$ and $q$ in $Q$ with a singularity (blowup) at the origin $z=0$, let us distinguish between two cases: $z=0$ is a Type I blowup if $g(v)<\infty$ and $z=0$ is Type II blowup if $g(v)=\infty$.

Following  paper \cite{Seregin2023}, we are going to study one of the possible  scenarios of Type II blowups generated  by two conditions. In the first one, it is supposed that there are  numbers $\varepsilon_0>0$ and $0<m_0<1$ and a sequence $r_k\downarrow0$ as $k\to\infty$ such that 
\begin{equation}
	\label{TypeII}
	M^{s,l}_{\kappa,m_0}(v,r_k)\geq \varepsilon_0
\end{equation}
for all natural $k$, where
$$M^{s,l}_{\kappa,m_0}(v,a)=a^{\kappa(1-m_0)}M^{s,l}_\kappa(v,a),
$$
$$M^{s,l}_\kappa(v,a)= \frac 1{a^\kappa}\int\limits_{-a^2} ^0\Big(\int\limits_{B(a)}|v|^sdx\Big)^\frac lsdt,$$
and numbers $s>1$ and $l>1$ satisfy restrictions:
\begin{equation}
	\label{restriction on s and l}
	l>\kappa:=l\Big(\frac 3s+\frac 2l-1\Big)>0.
\end{equation}
The second assumption in our scenario of Type II blowup is as follows: 
\begin{equation}
		\label{M1}
		M_1=\sup\limits_{0<r<1}\{A_{m_1}(v,r)+E_m(v,r)+D_m(q,r)\}<\infty
	\end{equation}
	for some $0<m<1$ with $m_1=2m-1$, where the quantities
$$E_m(v,r)=\frac 1{r^{m}}\int\limits_{Q(r)}|\nabla v|^2dz, \quad A_{m_1}(v,r)=\sup\limits_{-r^2<t<0}\frac 1{r^{m_1}}\int\limits_{B(r)}|v(x,t)|^2dx,$$
$$D_m(q,r)=\frac 1{r^{2m}}\int\limits_{Q(r)}|q|^\frac 32dz
$$
have been involved.

For more details about Type I and Type II blowups, we refer the reader to  papers \cite{Seregin2010}, \cite{Seregin2015}, and \cite{Seregin2023}.

The interesting question is about a relationship between numbers $m$ and $m_0$ that does not exclude  Type II blowup scenario described by \eqref{TypeII} and \eqref{M1}. In the second section of the note, see Theorem \ref{LPS1} there, we find a region for $m$ and $m_0$ that completely excludes scenario \eqref{TypeII} and \eqref{M1}.

Section 3 deals with a modification of the scenario \eqref{TypeII} and \eqref{M1}, where restriction \eqref{TypeII} is replaced with the new one that is   \eqref{one more}. The necessary condition for new scenario \eqref{one more} and \eqref{M1} to be happened is the existence of a non-trivial ancient solution to the Euler system in a certain class motivated by a particular scenario of the blowup.

In Section 4, several results on Liouville type theorems for ancient solutions to the Euler equations under the assumptions from the previous section are stated. The corresponding proofs can be found in papers \cite{Seregin2023} and \cite{Seregin2023-1}.

The final section devoted to the further analysis of the new scenario of Type II blowup with additional  Ladyzhenskaya-Prodi-Serrin type condition, see \eqref{LPS}, via Liouville type theorem. The reason, why \eqref{LPS} is marked by such a name, is as follows. If $m=1$, then restriction \eqref{M1} describes a potential Type I blowup and \eqref{LPS} becomes the Laddyzheskaya-Prodi-Serrin condition.

\setcounter{equation}{0}

\section{Impossible Cases of Potential Type II\\ Blowups in Scenario \eqref{TypeII} and \eqref{M1}    }

To state the main result of the section, let us introduce auxiliary numbers
\begin{equation}
	\label{numbers}
\frac{1}{	p(\eta)}=\frac \eta 6+\frac {3(1-\eta)} {10}, \quad \frac 1{q(\eta)}=\frac \eta 2+\frac {3(1-\eta)} {10},
\end{equation}
where $0\leq \eta\leq 1$ is a parameter. It is easy to check 
that 
$$1>\frac 3{p(\eta)}+\frac 2{q(\eta)}-1=\frac 12>0$$
for all $0\leq \eta \leq 1$.

\begin{theorem}
	\label{LPS1}
	Let $v$ and $q$ be a suitable weak solution to the Navier-Stokes equations in $Q$. It is supposed that the pair obeys  condition \eqref{M1}.

Assume further that, for some numbers $s$ and $l$, satisfying restriction (\ref{restriction on s and l}), the additional inequalities
\begin{equation}
	\label{numberscond}s <p(\eta),\quad l<q(\eta)\end{equation}
hold with some parameter $0\leq\eta \leq 1$.

Suppose that
\begin{equation}
\label{m_0}	
	0<m_0< \frac {\frac{3(\alpha-1)}2\Big(1-\frac s{p(\eta)}\Big)-\Big(\alpha l-\frac {3l}s-\alpha-1\Big)}{\frac{\alpha+1}2\kappa}=: f(m),
	\end{equation} 
	where $\alpha=2-m$.
	
	 Then 
\begin{equation}
	\label{no-typeII}
	\lim\limits_{r\to 0}M^{s,l}_{\kappa,m_0}(v,r)=0.
\end{equation}	
	\end{theorem}

\begin{remark} In paper \cite{Seregin2023},
it had been supposed that numbers $m$ and $m_0$ are related through identities:
\begin{equation}
	\label{Seregin2023}
	m=2-\alpha, \quad\alpha=\frac {l\Big(\frac 3s+1\Big)-(m_0-1)\kappa}{l\Big(\frac 3s+1\Big)+(m_0-1)\kappa} 
\end{equation}
Theorem \ref{LPS1} shows that restriction \eqref{Seregin2023} is too strong and excludes 
scenario of Type II described  by \eqref{TypeII} and \eqref{M1}. 
\end{remark}	
	Indeed,  finding $m_0(m)$ from 
	identities \eqref{Seregin2023}, one can easily verify that $m_0(m)<f(m)$ for any $0<m<1$. The latter means that $\varepsilon_0$ in growth condition \eqref{TypeII} must be equal to zero.

Now, let us start with our proof of Theorem \ref{LPS1}.
\begin{proof}
Suppose that statement \eqref{no-typeII} is wrong. Then, there are a positive number $\varepsilon_0$ and a sequence $r_k\to 0$ as $k\to\infty$ such that
\begin{equation}
	\label{non-statementLPS}
	M^{s,l}_{\kappa,m_0}(v,r_k)\geq \varepsilon_0>0
\end{equation}	
for all natural numbers $k$. 

Now, given $0<m<1$, define  the following scaling  
\begin{equation}
	\label{Euler-scaling}
	v^{\lambda,\alpha}(y,\tau)=\lambda^\alpha v(x,t)\qquad q^{\lambda,\alpha}(y,\tau)=\lambda^{2\alpha} q(x,t)\end{equation}
where $\alpha =2-m$ and 
$$x=\lambda y\qquad t=\lambda^{\alpha+1}\tau.$$
The choice of $\lambda$ is as follows:
$$\lambda=\lambda_k=r_k^\frac 2{\alpha+1}.$$
Then, after the change of variables, we arrive at the important inequality: 
\begin{equation}
	\label{non-trivial-2}
 \int\limits_{-1}^0\Big(\int\limits_{B(\lambda^{\frac {\alpha-1}2})}|v^{\lambda,\alpha}|^sdy\Big)^\frac lsd\tau=\lambda^\gamma M^{s,l}_{\kappa,m_0}(v,r_k)\geq \lambda^\gamma\varepsilon_0,
	\end{equation}
	where 
	$\gamma=\alpha l-\frac {3l}s-\alpha-1 +\frac {\alpha+1}2\kappa m_0.$

Next, setting $r=a\lambda $  for $a<1/\lambda$, we find, by the change of variables, the following estimate:
$$M_1\geq E_m(v,r)=\frac 1{r^m}\int\limits^0_{-r^2}\int\limits_{B(r)}|\nabla v|^2dxdt=$$$$=\frac {\lambda^{2-\alpha}}{r^m}\int\limits^0_{-(\frac r\lambda)^2/\lambda^{\alpha-1}}\int\limits_{B(\frac r\lambda)}|\nabla v^{\lambda,\alpha}|^2dyds=\frac 1{a^m}\int\limits^0_{-a^2\lambda^{1-\alpha}}\int\limits_{B(2a)}|\nabla v^{\lambda,\alpha}|^2dyds\geq
$$
$$\geq E_m(v^{\lambda,\alpha},a).$$ Arguing 
further in the same way, we have
\begin{equation}
	\label{M1scaling}
	M_1\geq \sup\limits_{0<a<1/\lambda}\{A_{m_1}(v^{\lambda,\alpha},a)+E_m(v^{\lambda,\alpha},a)+D_m(q^{\lambda,\alpha},a)\}.\end{equation}
Known multiplicative inequalities and arguments of the paper \cite{Seregin2023} lead to the existence of subsequences
of $v^{\lambda_k,\alpha}$ and $q^{\lambda_k,\alpha}$ such that:
\begin{itemize}
	\item $v^{\lambda_k,\alpha}\to u$ in $L_{3\nu}(Q(a))$;
	\item  $v^{\lambda_k,\alpha} {\stackrel{*}\rightharpoonup} u$ in $L_{2,\infty}(Q(a))$	\item $\nabla v^{\lambda_k,\alpha}\rightharpoonup \nabla u$ in $L_2(Q(a))$;
	\end{itemize}
for all $a>0$ and for all $1\leq \nu<\frac {10}9$. Moreover, the limit functions $u$ and $p$ possess the properties listed below:
$$	M_1\geq \sup\limits_{a>0}\{A_{m_1}(u,a)+E_m(u,a)+D_m(p,a)\}$$
	and
$$\partial_tu+u\cdot\nabla u=-\nabla p,\qquad {\rm div}\,u=0$$
in $Q_-$.

Now, we let  
$$s_0=6-\varepsilon,\quad s_1=\frac {10}3-\varepsilon.$$
For sufficiently small positive $\varepsilon$,  
$$\frac 1s> \frac 1{\tilde p(\eta)}=\frac \eta {s_0}+\frac {1-\eta}{s_1}>\frac 1{p(\eta)},\quad \frac 1l>\frac1{\tilde q(\eta)}=\frac \eta 2+\frac {1-\eta}{s_1}>\frac 1{q(\eta)}.
 $$
 Our next step is an application of  restriction \eqref{m_0}. Indeed, we can make a number $\varepsilon$ so small that  inequalities $$m_0<\bar m_0=\frac {\frac{3(\alpha-1)}2\Big(1-\frac s{\tilde p(\eta)}\Big)\frac ls-\Big(\alpha l-\frac {3l}s-\alpha-1\Big)}{\frac{\alpha+1}2\kappa}<f(m) $$ hold.
Repeating arguments of paper \cite{Seregin2023} and applying H\"older inequality, we obtain the bound $$\|v^{\lambda_k,\alpha}\|_{\tilde p(\eta),\tilde q(\eta),Q}\leq \|v^{\lambda_k,\alpha}\|_{s_0,2,Q}^\eta\|v^{\lambda_k,\alpha}\|_{s_1,Q}^ {1-\eta}\leq c_0<\infty$$ being valid 
all $k$.

 Making use of  H\"older inequality on more time, we find from \eqref{non-trivial-2} that:
$$\varepsilon_0\lambda^\gamma\leq c\lambda^{3\frac{\alpha-1}2\Big(1-\frac s{\tilde p(\eta)}\Big)\frac ls}\|v^{\lambda,\alpha}\|^l_{\tilde p(\eta),\tilde q(\eta),B(\lambda^{\frac{\alpha-1}2})\times]-1,0[}\leq c\lambda^{3\frac{\alpha-1}2\Big(1-\frac s{\tilde p(\eta)}\Big)\frac ls}c_0.$$
It remains to take into account definitions of $\gamma$ and $\bar m_0$ and get
$$\varepsilon_0\lambda^{\frac {\alpha+1}2\kappa(m_0-\bar m_0)}\leq c c_0.$$
Passing to the limit as $k\to\infty$, one can conclude that number $\varepsilon_0$ must vanish. This is a contradiction.
\end{proof}

\setcounter{equation}{0}

\section{One More Scenario of Type II Blowup}
Here, we replace assumption \eqref{TypeII}  by the following one: there are a positive number $\varepsilon_0$ and a sequence of $r_k\downarrow 0$ as $k\to\infty$ such that 
\begin{equation}	\label{one more}\overline{M}^{s,l}_{\kappa,m_0}(v,r_k):=r_k^{-\kappa m_0}\int\limits^0_{-r_k^{3-m}}\Big(\int\limits_{B(r_k)}|v|^sdx\Big)^\frac lsdt\geq \varepsilon_0>0
\end{equation}
for some $0<m<1$ and 

$$m_0=m_0(m)=\frac 1\kappa\Big(\alpha+1+\frac {3l}s-l\alpha\Big)$$
with $\alpha=2-m$. It is not so difficult  to check that $m_0(m)<1$ as long as $l>1$.

Obviously, since
\begin{equation}
	\label{imply Type II}
M^{s,l}_{\kappa,m_0}(v,r_k)\geq \overline{M}^{s,l}_{\kappa,m_0}(v,r_k),\end{equation}
condition \eqref{one more} describes  a potential Type II blowup. 
Necessary condition of the fact that such a  Type II blowup might happen follows from \eqref{m_0}, see Theorem \ref{LPS1}. Indeed, number $m\in ]0,1[$ must satisfy the inequality
\begin{equation}
	\label{m}
	m_0(m)\geq f(m).
\end{equation}

The latter is equivalent to the following one:
\begin{equation}
\label{F(m)}	F(m)=(l-1)m+3-2l+\frac {3l}{p(\eta)}\geq 0,\qquad m\in]0,1[.
\end{equation}
Simple calculations show  that there is a real number $m_*
$,  which is a unique root of the equation $F(m)=0$, and, moreover, 
$$ m_\star<1.$$ 
So, \eqref{F(m)} is equivalent to
\begin{equation}
\label{restriction-on-m}
m\in \mathcal M:=\{x\in \mathbb R: 
\,m_*\leq x\}\cap ]0,1[.
\end{equation}

Using the scaling introduced in the proof of Theorem \ref{LPS1}, we find the important identity
\begin{equation}
	\label{Non-trivial-identity}
	r_k^{-\kappa m_0}\int\limits^0_{-r_k^{3-m}}\Big(\int\limits_{B(r_k)}|v|^sdx\Big)^\frac lsdt=
	\int^0_{-1}\Big(\int\limits_B|v^{\lambda,\alpha}|^sdy\Big)^\frac ls d\tau=$$$$ =M^{s,l}_{\kappa,m_0}(v^{\lambda,\alpha},1)	\geq \varepsilon_0>0	
\end{equation}
with $\lambda=\lambda_k=r_k$. Then, repeating the same arguments as in the proof of Proposition 1.2 in \cite{Seregin2023}, we come to the  theorem.
\begin{theorem}
	\label{euler-non-trivial}
	Suppose that a pair $v$ and $q$ is a suitable weak solution to the Navier-Stokes equations in the unit space-time cylinder $Q$. Assume $v$ and $q$ satisfy the conditions \eqref{M1}, \eqref{one more}, and \eqref{restriction-on-m}.

Then, there are two functions $u$ and $p$ defined in $Q_-=\mathbb R^3\times]-\infty,0[$, with the following properties:
	$$\sup\limits_{a>0}\Big[\sup\limits_{-a^2<\tau<0}\frac 1{a^{m_1}}\int\limits_{B(a)}|u(y,\tau)|^2dy+\frac 1{a^{2m}}\int\limits_{Q(a)}|p|^\frac 32dyd\tau+$$
 \begin{equation}
	\label{basicestimates}
+\frac 1{a^{m}}\int\limits_{Q(a)}|\nabla u|^2dyd\tau\Big]\leq c<\infty;\end{equation}
\begin{equation}
	\label{Euler}\partial_\tau u+u\cdot\nabla u+\nabla p=0, \quad{\rm div}\,u=0
\end{equation}in $Q_-=\mathbb R^3\times ]-\infty,0[$ in the sense of distributions;
 
 for a.a. $\tau_0\in ]-\infty,0[$, the local energy inequality
$$\int\limits_{\mathbb R^3}|u(y,\tau_0)|^2\varphi(y,\tau_0)dy\leq$$
\begin{equation}
	\label{enerylocal}
\leq \int\limits_{-\infty}^{\tau_0}\int\limits_{\mathbb R^3} \Big(|u|^2\partial_\tau\varphi+u\cdot\nabla\varphi(|u|^2+2
p)
\Big)dyd
\tau
\end{equation}
holds for non-negative $\varphi\in C^\infty_0(\mathbb R^3\times \mathbb R)$;
 
 the function $u$ is not trivial in the sense
 \begin{equation}
	\label{nontrivial}
M^{s,l}_{\kappa,m_0}(u,1)
\geq \varepsilon_0/2.	
\end{equation} 
\end{theorem}


	
As it was shown in \cite{Seregin2023}, if 
$m<1/2(\Leftrightarrow m_1<0)$ then $u\equiv 0$ which contradicts with \eqref{nontrivial}. So, from now on, we may assume that 
$1>m\geq \max\{m_*,1/2\}$.


\setcounter{equation}{0}

\section{Some Liouville Type Theorem for the Euler system}	

Here, we are going to state several Liouville type theorems for the Euler system \eqref{Euler} whose ancient solutions satisfy local energy inequality \eqref{enerylocal} and bound \eqref{basicestimates}.

The simplest cases are:
\begin{itemize}
	\item $m<1/2$;
	\item ${\rm rot}\,{u}=0$. 
\end{itemize}

If any of the above conditions are met then $u=0$.
An elementary proof is given  in \cite{Seregin2023}. 

One can look for non-trivial solutions to Euler system \eqref{Euler} in self-similar form. To this end, let us introduce profile functions $U$ and $P$ so that:
$$u(x,t)=\frac 1{(-t)^{\frac \alpha{\alpha+1}}}U(y,\tau),\quad p(x,t)=\frac 1{(-t)^{\frac {2\alpha}{\alpha+1}}}P(y,\tau),
$$
where 
$$y=\frac x{(-t)^{\frac 1{\alpha+1}}}, \quad \tau =-\ln(-t)
$$
and the number $\alpha=2-m>1$. 
 Then, the Euler system takes the form 
\begin{equation}
	\label{special-form}
	\partial_\tau U+\frac \alpha{\alpha+1}U+\frac 1{\alpha+ 1}y\cdot\nabla U+U\cdot\nabla U+\nabla P=0,\quad{\rm div}\,U=0
\end{equation}
for  $(y,\tau)\in\mathbb R^3\times \mathbb R$.	

Consider the case  $U(y,\tau)=U(y)$, i.e., $u$ is a really self-similar solution.  Condition \eqref{basicestimates} gives a certain restriction on the function class for the profile function $U$:
$$
\sup\limits_{b\geq 1}\frac 1{b^{m_1}}\int\limits_{B(b)}|U(y)|^2dy+\sup\limits_{b>0}\frac 1{b^m}\int\limits_{B(b)}|\nabla U|^2dy+$$
\begin{equation}
	\label{1sts-s}
	+\sup\limits_{b>0}\frac 1{b^{2m}}\int\limits_{B(b)}|P|^\frac 32dy\leq c<\infty.
\end{equation}

\begin{pro}
	\label{simpleLiouville}  
	Suppose that the following additional conditions hold:
	\begin{equation}
		\label{1staddition}
	\frac 12< m	<\frac 35\quad(\Leftrightarrow 0<m_1< \frac 15)
	\end{equation}
	and
	\begin{equation}
		\label{2ndaddition}
		\sup\limits_{b\geq 1}\frac 1{b^{\gamma m_1}}\int\limits_{B(b)}|U|^2dy<\infty
	\end{equation}
	with $$0\leq \gamma < \frac {\ln {(2+m_1)}-\ln2}{m_1\ln 2}.$$
	Then $U(y)\equiv0$.
\end{pro}
Similar arguments work for discrete self-similar solutions.


Axisymmetric solutions are quite interesting because it is known that they do not have Type I singularities, see \cite{Seregin2020}.

To study axisymmetric case, we are going to exploit  cylindrical coordinates $r=|x'|$, with $x'=(x_1,x_2,0)$, $\vartheta$, and $x_3$. So, $u=u_re_r+u_\vartheta e_\theta+u_3e_3$, where $e_r$, $e_\vartheta$, and $e_3$ form the orthonormal  basis of the cylindrical coordinates.

In what follows, just for convenience, we are going to replace  all spatial balls $B(a)$ with spatial cylinders $\mathcal C(a):=\{|y'|<a,\, |y_3|<a\}$.

As it has been shown in \cite{Seregin2023-1}, the limit Euler system must have zero swirl, i.e., $u_\theta=0$, and thus the Euler system is reduced to following one:
$$\partial_tu_r+u_ru_{r,r}+u_3u_{r,3}+p_{,r}=0,$$
\begin{equation}
	\label{velocity}
	\partial_tu_3+u_ru_{3,r}+u_3u_{3,3}+p_{,3}=0,
\end{equation}
$$\frac 1r(ru_r)_{,r}+u_{3,3}=0$$
or equivalently 
$$\partial_tf+u_rf_{,r}+u_3f_{,3}=0,$$
\begin{equation}
	\label{vorticity}
	\Delta \psi-\frac 2r\psi_{,r}=r^2f,
\end{equation}
$$u_r=\frac 1r\psi_{,3},\qquad u_3=-\frac 1r\psi_{,r}.$$
It is worth to notice that the only one component of the vorticity 
$$\omega_\vartheta(u)=u_{r,3}-u_{3,r} =rf$$
is not necessary to be vanishing, i.e., 
$\omega(v)=\omega_\vartheta(v)e_\vartheta$.

Assume, in addition, the following two condition hold. The first one is:
 \begin{equation}
 	\label{secondderiveNsl} 
 	N^{s_1,l_1}(u,a):=\frac 1{a^{\gamma_*l_1}}\int\limits^0_{-a^2}\Big(\int\limits_{B(a)}(|\nabla^2u|^{s_1}+|\partial_tu|^{s_1})dx\Big)^\frac {l_1}{s_1}dt\leq c<\infty.
 	\end{equation}
It is valid for all $a>0$. Numbers $s_1,l_1>1$ obey the identity
\begin{equation}
	\label{for second deriv}
	\frac 3{s_1}+\frac 2{l_1}=4
\end{equation}
and 
$\gamma_*(s_1,m)=1-(3/(2s_1)-1)(1-m).
 $

In the second one, velocity field $u$ is supposed to be bounded outside of the axis of
symmetry, i.e., 
	\begin{equation}
		\label{decaycond1}
	|u(x,t)|\leq \frac c{|x'|^\alpha}	\end{equation}
	for all $z=(x,t)\in Q_-
	$ with $
	|x'|>0$.    Here, as usual, $\alpha=2-m$. 
 
 Why the above restrictions on the function class of the Liouville type theorem are reasonable is explained in the preprint \cite{Seregin2023-1}.
 
 The next statement is a kind of conservation law for ancient solutions to the Euler system.

\begin{pro} Assume that a pair $u$ and $p$ satisfies 
the Euler equations \eqref{Euler},  local energy inequality \eqref{enerylocal} and bound \eqref{basicestimates}. 
Suppose further that conditions \eqref{secondderiveNsl}	and \eqref{decaycond1} are met as well and 
	numbers $s_1$ and $l_1$ obey additional restrictions
$$l_1\leq s_1,\qquad 	m<\frac {4l_1-3}{l_1+1}.$$

Assume, further, that there is a number $t_0\leq 0$ such that 
	\begin{equation}
		\label{assumption}
		g(t_0):=\frac 2{l_1}\int\limits_{\mathbb R^3}\Big(\frac {|\omega_\vartheta(u(x,t_0)|}{r}\Big)^\frac {l_1}2dx<\infty.
	\end{equation}

Then 
\begin{equation}
	\label{preservation}
	g(t):=\frac 2{l_1}\int\limits_{\mathbb R^3}\Big(\frac {|\omega_\vartheta(u(x,t)|}{r}\Big)^\frac {l_1}2dx=g(t_0)
\end{equation}	
for all $t\leq 0$.

In particular, if $\omega_\theta(\cdot,t_0)=0$ in $\mathbb R^3$, then $u\equiv0$ in $Q_-$.
\end{pro}
For the detail proof, we refer the reader to preprint \cite{Seregin2023-1}.


\setcounter{equation}{0}

\section{Ladyzhenskaya-Prodi-Serrin Type Condition}

Let us describe  standing assumptions of this section:
\begin{itemize}
	\item $v\in L_3(Q)$ and $q\in L_\frac 32(Q)$ satisfy the Navier-Stokes equations $$
	\partial_tv+v\cdot\nabla v-\Delta v=-\nabla q,\qquad {\rm div}\,v=0$$ in $Q$ in the sense of distributions;
	\item $v\in L_\infty(B\times ]-1,-a^2[\cup \{r_1<|x|<1\}\times ]-1,0[)$ for any $a\in ]0,1[$ and for some $r_1\in ]0,1[$.
\end{itemize}
Typically, these conditions appear in a scenario  of the first time blowup.  It is easy to show 
that the pair $v$ and $q$ is a suitable weak solution to the Navier-Stokes equations in $Q(b)$ for any $0<b<1$, see for example \cite{SS2009}. So, without loss of generality, we may assume that  the pair $v$ and $q$ is a suitable weak solution in $Q$.

In this section, we are going to prove the following statement.
\begin{theorem}
	\label{LPS2}
	Let a pair $v$ and $q$ satisfy the above standing assumptions. In addition, it is supposed that it obeys the following conditions:

(i) bound	 \eqref{M1} 
	for some $\frac 12\leq m<1$ with $m_1=2m-1$;

(ii)	additional restriction
\begin{equation}
\label{M2}	
	M_2=\sup\limits_{0<r<1}N^{s_1,l_1}(v,r)<\infty
	\end{equation} for some $1<s_1,l_1<\infty$ such that $3/s_1+2/l_1=4$, for definition of the quantity $N^{s,l}$ see \eqref{secondderiveNsl};

(iii) Ladyzhenskaya-Prodi-Serrin type assumption
\begin{equation}
	\label{main LPS condition}
	v\in L_{s_2,l_2}(Q)
\end{equation} for some numbers $1\leq s_2\leq \infty$ and $1\leq l_2\leq \infty$ subjected to the identity
\begin{equation}
	\label{restrictionLPS}
	\frac 3{s_2}+\frac {\alpha+1}{l_2}= \alpha
\end{equation}	 with $\alpha=2-m$.
If $m=\frac 12\, (\Leftrightarrow m_1=0\Leftrightarrow\alpha=\frac 32)$ and $l_2=\infty$, assume, additionally, 
\begin{equation}
	\label{continuity}
	v\in C([-\delta^2_1,0];L_2(B(\delta_1)))
\end{equation}
for some number $0<\delta_1\leq1$.

Then 
\begin{equation}
	\label{statementLPS}
	\lim\limits_{r\to 0}\overline M^{3,3}_{2,m}(v,r)=0.
\end{equation}	
	\end{theorem}

\begin{remark}
	\label{LPS} The above theorem states that the scenario, described in Section 3 for a particular case $s=l=3$, $\kappa=2$, is impossible under the additional assumption \eqref{LPS}. We call it a Ladyzhenskaya-Prodi-Serrin type condition since, for $m=1(\Leftrightarrow \alpha=1)$, identity \eqref{restrictionLPS} becomes exactly the celebrated Ladyzhenskaya-Prodi-Serrin condition.
\end{remark}
\begin{proof}
Suppose that statement \eqref{statementLPS} is wrong. Then there are a positive number $\varepsilon_0$ and a sequence $r_k\downarrow 0$ as $k\to  \infty$ such that
\begin{equation}
	\label{non-statementLPS}
	\overline M^{3,3}_{2,m_0}(v,r_k)\geq \varepsilon_0
\end{equation}	
for all natural numbers $k$.

We can repeat all the arguments of the proof of Theorem \ref{one more} based on the Euler scaling
\eqref{Euler-scaling} for a particular case $s=l=3$, $\kappa=2$, $\lambda=\lambda_k=r_k$.

Then, after the change of variables, we arrive at the identity: 
$$\overline M^{3,3}_{2,m_0}(v,r_k)= \int\limits_{Q}|v^{\lambda,\alpha}|^3dyds.$$
The latter leads to the following version of condition \eqref{non-statementLPS}:
\begin{equation}
	\label{non-trivial2}
0<	\varepsilon_0\leq \int\limits_{Q}|v^{\lambda,\alpha}|^3dyds.
	\end{equation}

By Theorem \ref{one more},  the following facts are valid (up to selecting a subsequence) : 
\begin{itemize}
	\item $v^{\lambda_k,\alpha}\to u$ in $L_{3\nu}(Q(a))$;
	\item  $v^{\lambda_k,\alpha} {\stackrel{ *}\rightharpoonup} u$ in $L_{2,\infty}(Q(a))$	\item $\nabla v^{\lambda_k,\alpha}\rightharpoonup \nabla u$ in $L_2(Q(a))$;
	\item  $q^{\lambda_k,\alpha}\rightharpoonup p$ in $L_{\frac 32}(Q(a))$;	
	\item $\nabla^2v^{\lambda_k,\alpha}\rightharpoonup\nabla^2u$ in $L_{s_1,l_1}(Q(a))$;
	\item $\partial_tv^{\lambda_k,\alpha}\rightharpoonup\partial_tu$ in $L_{s_1,l_1}(Q(a))$	\end{itemize}
for all $a>0$ and for all $1\leq \nu< \frac {10}3$. Moreover, the limit functions $u$ and $p$ possess the properties:
\begin{equation}
	\label{M1scalinglimit}
	M_1\geq \sup\limits_{0<a}\{A_{m_1}(u,a)+E_m(u,a)+D_m(p,a)\};\end{equation} 
\begin{equation}
	\label{M2scalinglimit}
M_2\geq\sup\limits_{0<a}
	N^{s_1,l_1}(u,a);
	\end{equation}
$$\partial_tu+u\cdot\nabla u=-\nabla p,\qquad {\rm div}\,u=0$$
in $Q_-$; 
\begin{equation}
	\label{non-trivial-limit}
	0<\varepsilon_0\leq \int\limits_Q|u|^3dxds.
\end{equation}

	Now, let us consider the easy case
	\begin{equation}
		\label{easy-case}
		l_2<\infty.
	\end{equation}
Here, we have
$$\int\limits_Q|v^{\lambda,\alpha}|dyds=
\frac 1{\lambda^4}\int\limits^0_{-\lambda^{\alpha+1}}\int\limits_{B(\lambda)}|v|dxdt\leq 
$$	
$$\leq \Big(\int\limits^0_{-\lambda^{\alpha+1}}\Big(\int\limits_{B(\lambda)}|v|^{s_2}dx\Big)^\frac {l_2}{s_2}dt\Big)^\frac 1{l_2}\to 0$$
as $\lambda\to 0$. Since $v^{\lambda,\alpha}\to u$ a.e. in $Q$, function $u$ must be zero in $Q$. So, we arrive at the contradiction with \eqref{non-trivial-limit}.

So, now, one should consider the case $l_2=\infty$ ($\Rightarrow s_2=3/\alpha\geq 2$) carefully. 	
	
First, by scaling, we find the estimate
$$\sup\limits_{-1<a<0}\int\limits_B|v(x,t)|^\frac 3\alpha dx=\sup\limits_{-1/\lambda^{\alpha+1}<s<0}\int\limits_{B(1/\lambda)}|v^{\lambda,\alpha}(y,s)|^\frac 3\alpha dy \geq
$$	
$$\geq \sup\limits_{-b^2<s<0}\int\limits_{B(a)}|v^{\lambda,\alpha}(y,s)|^\frac 3\alpha dy
$$
being valid for any positive numbers $a$ and $b$ and sufficiently small $\lambda$. It yields:
\begin{equation}
	\label{additional-est}
	\|u\|_{\frac 3\alpha,\infty,Q_-}\leq 
	\|v\|_{\frac 3\alpha,\infty,Q}.
\end{equation}	

Now, let us list some useful properties of $u$ coming from inequality \eqref{additional-est} and estimate \eqref{M2scalinglimit}:
\begin{itemize}
	\item function $t\mapsto\int\limits_{B(a)}u(x,t)\cdot w(x)dx$ is continuous on $[-a^2,0]$ for each $w\in L_{\frac 3{3-\alpha}}(B(a))$;
	\item $\|u(\cdot,t)\|_{\frac 3\alpha,B(a)}<\infty$ for each $t\in [-a^2,0]$;
	\item $u\in C([-a^2,0];L_s(B(a))$ for any $s<3/\alpha$
\end{itemize}
for any $a>0$.

For the same reason, 
we may assume also that:
\begin{equation}
	\label{est-for-v}
	\|v(\cdot,t)\|_{\frac 3\alpha,B(3/4)}<\infty
\end{equation}	
for all $
-(3/4)^2\leq t\leq 0$.	

Next, for any test function $\varphi\in C^\infty_0(B(r))$,
$$\Big|\int\limits_{B(r)}v^{\lambda,\alpha}(y,0)\cdot\varphi(y)dy\Big|=\Big|\lambda^{\alpha-3}\int\limits_{B(r\lambda)}v(x,0)\cdot\varphi(x/\lambda)dx\Big|\leq 
$$	
$$\leq	c\|\varphi\|_{\infty,\mathbb R^3}
\Big(\int\limits_{B(r\lambda)}|v(x,0)|^\frac 3\alpha dx\Big)^\frac \alpha 3\to 0
$$	provided $\lambda\to0$. Standard arguments allow us to conclude that, for all $a>0$,
\begin{equation}
	\label{weak-continuous-conv}
	\int\limits_{B(a)}v^{\lambda,\alpha}(y,
	\tau)\cdot\varphi(y)dy\to \int\limits_{B(a)}u(y,
	\tau)\cdot\varphi(y)dy
	\end{equation}
in $C([-a^2,0])$ for all $\varphi\in L_{s'_1}(B(a))$, where $s'_1=s_1/(s_1-1)$. The latter implies that
\begin{equation}
	\label{last-moment-zero}
	u(\cdot,0)=0
\end{equation}	
in $\mathbb R^3$.

Now, our goal is to establish a certain local energy 	identity. To this end, we pick up a test function $\psi\in C^\infty_0(Q_-)$. If let $\varphi_\lambda(x,t)=\psi(x/\lambda, t/\lambda^{\alpha+1})$, then $\varphi_\lambda \in C^\infty_0(Q)$ provided $\lambda a<1$, where a number $a$ is chosen so large that 
${\rm spt}\,\psi\subset Q(a)$. Then, according to the standing assumptions,
we have 
$$\int\limits_Q\Big(-|v|^2(\partial_t\varphi_\lambda+\Delta\varphi_\lambda)+2\varphi_\lambda|\nabla v|^2-(|v|^2+2q)v\cdot\nabla\varphi_\lambda\Big)dxdt=0.
$$
Scaling 	in the latter identity gives us another one:
	$$\int\limits_{Q(a)}\Big(-|v^{\lambda,\alpha}|^2(\partial_s\psi+\lambda^{\alpha-1}\Delta\psi)+\lambda^{\alpha-1}2\psi|\nabla v^{\lambda,\alpha}|^2-(|v^{\lambda,\alpha}|^2+$$$$+2q^{\lambda,\alpha})v^{\lambda,\alpha}\cdot\nabla\psi\Big)dyds=0.
$$
Here, we can take the limit passing $\lambda\to0$ and get the following 
local energy identity for the Euler equations:
\begin{equation}
	\label{EulerIdenity}
\int\limits_{Q_-}\Big(|u|^2\partial_t\psi+(|u|^2+2p)u\cdot\nabla \psi\Big)dyds=0
\end{equation}	
 for any $\psi\in C^\infty_0(Q_-)$.


Our aim now to estimate the second and third terms \eqref{EulerIdenity}. To this end,  given $a>0$, let us start
 with the cubic term. By 	
	  interpolation, Gagliardo-Nirenberg inequality, and \eqref{additional-est},  we show the following:
$$\int\limits_{B(a)}|u|^3dx\leq \Big(\int\limits_{B(a)}|u|^\frac 3\alpha dx\Big)^\frac {\alpha}{2\alpha-1}\Big(\int\limits_{B(a)}|u|^6dx)^\frac {\alpha-1}{2\alpha-1}\leq$$
$$\leq c(v)\Big(\int\limits_{B(a)}|\nabla u|^2dx+\frac 1{a^2}\int\limits_{B(a)}|u|^2dx\Big)^{\frac 12 6\frac {\alpha-1}{2\alpha-1}} \leq $$
$$\leq  c(v)\Big(\int\limits_{B(a)}|\nabla u|^2dx+\frac {a^{m_1}}{a^2}A_{m_1}(u,a)\Big)^\frac {3(\alpha-1)}{2\alpha-1}.$$
With the help of \eqref{M1scalinglimit},  for $a^2>\max\{1,T\}$, it is possible to find 
$$\frac 1a\int\limits^0_{-T}\int\limits_{B(2a)}|u|^3dx dt\leq 
\frac {c(v)}a\int\limits^0_{-T}\Big(\int\limits_{B(2a)}|\nabla u|^2dx+\frac {(2a)^{m_1}}{(2a)^{2}}A_{m_1}(u,2a)\Big)^\frac {3(\alpha-1)}{2\alpha-1}dt\leq $$
$$\leq \frac {c(v)}a\Big[(2a)^mE_m(u,2a)+(2a)^{m_1}A_{m_1}(u,2a)\Big]^\frac {3(\alpha-1)}{2\alpha-1}T^{1-\frac {3(\alpha-1)}{2\alpha-1}}\leq 
$$
$$\leq c(v,\alpha,T)a^{m\frac {3(\alpha-1)}{2\alpha-1}-1}M_{1}^\frac {3(\alpha-1)}{2\alpha-1}.$$
It is easy to verify that $m\frac {3(\alpha-1)}{2\alpha-1}=\frac {3(2-\alpha)(\alpha-1)}{2\alpha-1}<1$ if $1<\alpha<2$ and thus 
\begin{equation}
	\label{est-for-u}\lim\limits_{a\to \infty}\frac 1a\int\limits_{-T}^0\int\limits_{B(a)}|u|^3dx dt=0\end{equation}
for each fixed $T>0$.

Now, our task is to evaluate the third term in the left hand side of \eqref{EulerIdenity}.
To this end, first we are going to use a particular consequence of \eqref{M1scalinglimit}:
\begin{equation}
	\label{particular}
	\sup\limits_{a>1}D_m(p,a)\leq M_{1}<\infty.
\end{equation}

Second, notice that $p$ is a solution to the pressure equation: 
$$-\Delta p={\rm div\,div}\,u\otimes u.$$
Usual way to treat the pressure $p$ is to split it into two parts:
$$p=p_1+p_2,
$$
where 
$$p_1(x,t)=-\frac 13|u(x,t)|^2+\frac 1{4\pi}\int\limits_{\mathbb R^3}K(x-y):u(y,t)\otimes u(y,t)dy,$$
$$ K(y)=\nabla^2\Big(\frac 1{|y|}\Big)$$
for $|y|>0$, so that,  by the standard singular integral estimate and by \eqref{additional-est},
$$\|p_1\|_{\frac 3{2\alpha},\infty, Q_-}\leq c\|v\|_{\frac 3{\alpha},\infty, Q},$$
and 
$$\Delta p_2(\cdot,t)=0$$
in $\mathbb R^3$.

In turn, for $x\in B(2a)$, the first term $p_1$ can be presented as sum  of three terms so that
$$p_1(x,t)=p_{11}(x,t)+p_{12}(x,t)+c_0(t),$$  
where
$$p_{11}(x,t)=-\frac 13|u(x,t)|^2+\frac 1{4\pi}\int\limits_{B(4a)}K(x-y):u(y,t)\otimes u(y,t)dy,$$
$$p_{12}(x,t)=\frac 1{4\pi}\int\limits_{\mathbb R^3\setminus B(4a)}(K(x-y)-K(0-y)):u(y,t)\otimes u(y,t)dy,$$
and 
$$c_0(t)=\frac 1{4\pi}\int\limits_{\mathbb R^3\setminus B(4a)}K(0-y):u(y,t)\otimes u(y,t)dy.$$

Let us start with evaluation of $c_0(t)$. Indeed, we have:
$$|c_0(t)|\leq c\int\limits_{\mathbb R^3\setminus B(4a)}\frac 1{|y|^3}|u(y,t)|^2dy\leq c \|u\|^2_{\frac 3\alpha,\infty,Q_-}\times$$
$$\times\Big(\int\limits_{\mathbb R^3\setminus B(4a)}\Big(\frac 1{|y|^3}\Big)^\frac 3{3-2\alpha}dy\Big)^\frac {3-2\alpha}3\leq$$
$$\leq c(\alpha) \|u\|^2_{\frac 3\alpha,\infty,Q_-}a^{-2\alpha}$$
for all $t\leq 0$. The above arguments are valid under the assumption that $\alpha<3/2$. If $\alpha=3/2 (\Leftrightarrow m=1/2)$, the proof is even easier.

Next, the term $p_{11}$ can be easily estimated with the help of singular integral theory:
$$\|p_{11}(\cdot,t)\|_{\frac 32,B(2a)}\leq c\|u\|^2_{3,B(4a)}.$$
As to $p_{12}$, we are going to use known arguments, see for example \cite{KikSer}, and show that
$$|p_{12}(x,t)|\leq c\sum\limits^\infty_{i=0}\frac 1{(2^{i+2}a)^4}\int\limits_{B(2^{i+3}a)}|u(y,t)|^2dy\leq $$
$$\leq c\sum\limits^\infty_{i=0}\frac 1{(2^{i+2}a)^4}\Big(\int\limits_{B(2^{i+3}a)}|u(y,t)|^\frac 3\alpha dy\Big)^\frac{2\alpha}3(2^{i+3}a)^{3(1-\frac {2\alpha}3)}\leq $$
$$\leq c(\alpha) \|u\|^2_{\frac 3\alpha,\infty,Q_-}a^{-1-2\alpha}$$
for all $t\leq 0$ and for all $x\in B(2a)$.
Therefore,
$$\|p_1(\cdot,t)\|_{\frac 32, B(2a)}^\frac 32\leq c\|u(\cdot,t)\|^3_{3,B(4a)}+c \|u\|^3_{\frac 3\alpha,\infty,Q_-}a^{-(1+2\alpha)\frac 32}a^3+
$$
$$+c \|u\|^3_{\frac 3\alpha,\infty,Q_-}a^{3-3\alpha}.$$
So, after integration in time, we find
$$\frac 1a\int\limits^0_{-T_0}\int\limits_{B(2a)}|p_1|^\frac 32de\leq $$
\begin{equation}
	\label{est-for-p1}
\leq 	\frac ca
\int\limits^0_{-T_0}
\int\limits_{B(4a)}|u|^3de+\|u\|^3_{3/\alpha,\infty,Q_-}T_0a^{3(1-\alpha)-1}(1+a^{-3/2})\to 0
\end{equation}
as $a\to \infty$.

For the second counter-part of the pressure, we have
$$\frac 1{a^{2m}}\int\limits^0_{-T_0}\int\limits_{B(2a)}|p_2|^\frac 32de\leq c\frac 1{a^{2m}}\int\limits^0_{-T_0}\int\limits_{B(2a)}|p|^\frac 32de+c\frac 1{a^{2m}}\int\limits^0_{-T_0}\int\limits_{B(2a)}|p_1|^\frac 32de\leq $$
$$
\leq cD_m(p,a)+\frac c{a}\int\limits^0_{-T_0}\int\limits_{B(2a)}|p_1|^\frac 32de
$$
provided $a^2\geq \max\{T_0,1\}$ and $2m \geq1$.
So, the latter gives:
\begin{equation}
	\label{est-for-p2}
C_*(v,p,T_0,m):=\sup\limits_{a^2\geq \max\{T_0,1\}}	\frac 1{a^{2m}}\int\limits^0_{-T_0}\int\limits_{B(2a)}|p_2|^\frac 32de<\infty.
	\end{equation}

Now, our task is to show that   estimate \eqref{est-for-p2}  implies $p_2\equiv 0$. To this end, fix numbers $-\infty<-T_0<a<b\leq 0$ and introduce a function 
$$g(x_0):=\int\limits^b_ap_2(x_0,t)dt.$$ This function is harmonic in $\mathbb R^3$ and  thus we have 
$$|g(x_0)|\leq \int\limits^b_a|p_2(x_0,t)| dt\leq \int\limits^b_a\frac 1{|B(|x_0|)|}\int\limits_{B(x_0,|x_0|)}|p_2(x,t)|dxdt \leq$$
$$\leq 2^3 \int\limits^b_a\frac 1{|B(2|x_0|)|}\int\limits_{B(2|x_0|)}|p_2(x,t)|dx dt\leq $$
$$\leq 2^3 \int\limits^b_a\Big(\frac 1{|B(2|x_0|)|}\int\limits_{B(2|x_0|)}|p_2(x,t)|^\frac 32dx \Big)^\frac 23dt\leq $$
$$\leq 2^3 |b-a|^\frac 23\Big(\frac 1{|B(2|x_0|)|}\int\limits^b_a\int\limits_{B(2|x_0|)}|p_2(x,t)|^\frac 32dx dt\Big)^\frac 23.$$
Assuming that $(2|x_0|)^2>T_0$, we find 
$$|g(x_0)|\leq c( C_{*})(b-a)^\frac 23\Big (\frac 1{|x_0|}\Big)^\frac {2(3-2m)}3\to 0
$$ as $|x_0|\to \infty$. So, $g(x_0)=0$ and, by arbitrariness of $a$ and $b$, we conclude that $p_2$ is identically equal to zero in $\mathbb R^3\times ]-T_0,0[$ for any $T_0>0$.

Evaluating the third term on the right hand side of \eqref{EulerIdenity}, one can show
$$\frac 1a\int\limits^0_{-T_0}\int\limits_{B(a)}|p||u|dyds\leq $$
$$\leq  c\Big(\frac 1a\int\limits^0_{-T_0}\int\limits_{B(2a)}|u|^3dx dt\Big)^\frac 13\Big(\frac 1a\int\limits^0_{-T_0}\int\limits_{B(2a)}|p_{1}|^\frac 32dx dt\Big)^\frac 23 \to0 $$
as $a\to\infty$, see \eqref{est-for-u} and   \eqref{est-for-p1}.

Now, consider the first sub-case when
$\alpha <3/2$ and thus $2< 3/\alpha$. Hence, $u\in C([-a^2,0];L_2(B(a))$ for any $a>0$ and it makes possible to deduce from \eqref{EulerIdenity} the identity:
$$\int\limits_{\mathbb R^3}|u(y,s_2)|^2\varphi_a(y)dy-\int\limits_{\mathbb R^3}|u(y,s_1)|^2\varphi_a(y)dy+$$
$$+\int\limits^{s_2}_{s_1}\int\limits_{\mathbb R^3}(|u|^2+2p)u\cdot\nabla\varphi_a	dy ds=0$$
for any $-a^2\leq s_1\leq s_2\leq 0$ and for a cut-off function $\varphi_a\in C^\infty_0(B(a))$ such that: $0\leq \varphi_a\leq 1$ in $B(a)$, $\varphi_a=1$ in $B(a/2)$, and $|\nabla \varphi_a|\leq c/a$ in $B(a)$.
	
Setting $s_2=0$ and $s_1=-T_0$, we find	
\begin{equation}
	\label{main-identity}
\int\limits_{\mathbb R^3}|u(y,-T_0)|^2\varphi_a(y)dy=\int\limits^{0}_{-T_0}\int\limits_{\mathbb R^3}(|u|^2+2p)u\cdot\nabla \varphi_a	dy d\tau	\end{equation}
for any $a>0$. 
As a result,  
$$\int\limits_{\mathbb R^3}|u(y,-T_0)|^2\varphi_a(y)dy\to0$$ as $a\to \infty$. From Beppo-Levi theorem, it follows that $u(\cdot,-T_0)=0$ in $\mathbb R^3$ for any $T_0>0$. The latter contradicts with \eqref{non-trivial-limit}.

Now, let us consider the second sub-case  when $l_2=\infty$, $\alpha=\frac{3}{2}$ which implies $m=1/2$, $m_1=0$, and $s_2=2$ 
and thus $u\in L_{2,\infty}(Q_-)$. By arbitrariness of $T_0$ and function $\chi$, one can deduce that 
\begin{equation}
	\label{constantC*}
g(s)=\int\limits_{\mathbb R^3}|u(y,s)|^2dy
=C_*=constant\end{equation}
for $s<0$. If  continuity of $g$ at $s=0$ would take a place, our proof would be finished. Here, our arguments are going   to exploit condition \eqref{continuity} in the following way. By that condition, given $\varepsilon>0$, there is a number $\tau_1>0$ such that
$$\|v(\cdot,t)-v(\cdot,0)\|_{2,B(\delta_1)}<\varepsilon
$$ as long as $|t|<\tau_1$. Using our scaling, we find
$$\int\limits_{B(\delta_1/\lambda)}|v^{\lambda,\alpha}(y,t/\lambda^{\alpha+1})-v^{\lambda,\alpha}(y,0)|^2dy<\varepsilon^2.
$$
Now, given $a>0$, it is possible to make $\lambda$ so small that
$$\frac {\delta_1}\lambda> a, \qquad \frac {\tau_1} {\lambda^{\alpha+1}}>a^2.
$$ Therefore, the following is true:
$$\int\limits_{B(a)}|v^{\lambda,\alpha}(y,s)-v^{\lambda,\alpha}(y,0)|^2dy<\varepsilon^2
$$
for any $-a^2\leq s\leq 0$. Taking an arbitrary function $\chi\in C^\infty_0([-a^2,0])$ and integrating the latter inequality in $s$, we derive the following one:
\begin{equation}
	\label{integral-inequality}
\int\limits_{Q(a)}\chi^2(s)|v^{\lambda,\alpha}(y,s)-
v^{\lambda,\alpha}(y,0)|^2dy ds\leq \varepsilon^2\int\limits^0_{-a^2}\chi^2(s)ds.
\end{equation} Next, we know 
$$v^{\lambda,\alpha}\to u
$$ in $L_2(Q(a))$ and
$$\int\limits_{Q(a)}|v^{\lambda,\alpha}(y,0)|^2dy ds=a^2\int\limits_{B(a)}|v^{\lambda,\alpha}(y,0)|^2dy=$$$$=a^2\int\limits_{B(\lambda a)}|v(y,0)|^2dy\leq a^2\|v\|_{2,\infty,Q}^2.
$$ So, we may think that $v^{\lambda,\alpha}(\cdot,0)$
converges weakly in $L_2(Q(a))$ and the limit is $u(\cdot,0)=0$, see  \eqref{weak-continuous-conv} and \eqref{last-moment-zero}.

Now, taking the limit in \eqref{integral-inequality} as $\lambda\to 0$, one can  conclude that 
$$\int\limits_{-T_0}^0\chi^2(s)\int\limits_{B(a)}|u(y,s)|^2dyds\leq \varepsilon^2 \int\limits_{-T_0}^0\chi^2(s)ds
$$ for $a^2>T_0$. Passing to the limit as $a\to \infty$ and using \eqref{constantC*}, we find
$$C_*\int\limits_{-T_0}^0\chi^2(s)ds= \int\limits_{-T_0}^0\chi^2(s)\int\limits_{\mathbb R^3}|u(y,s)|^2dyds\leq \varepsilon^2 \int\limits_{-T_0}^0\chi^2(s)ds.
$$ By arbitrariness of $\varepsilon$, $C_*=0$.
\end{proof}

\end{document}